\documentclass[11pt,reqno]{amsart}

\usepackage{tikz}
\usetikzlibrary{decorations.pathreplacing}
\usepackage{psfrag}
\psfrag{A}{\fontsize{24}{24}$\psi_0^\ve$}
 \psfrag{B}{\fontsize{24}{24}$\ov{\psi_0^\ve}$}

\psfrag{1}{\fontsize{16}{16}$k_n$}
\psfrag{2}{\fontsize{16}{16}$k_{n-1}$}
\psfrag{3}{\fontsize{16}{16}$k_{n-2}$}
\psfrag{n}{\fontsize{16}{16}$k_1$}
\psfrag{o}{\fontsize{24}{24}$\mathcal{O}_\varepsilon$}
\psfrag{j}{\fontsize{20}{20}$t$}
\psfrag{q}{\fontsize{16}{16}$v + \xi / 2$}
\psfrag{w}{\fontsize{16}{16}$v - \xi / 2$}
\psfrag{r}{\fontsize{16}{16}$v - \xi / 2 + \sum_{l=1}^n k_l$}
\psfrag{0}{\fontsize{16}{16}$0$}
\psfrag{e}{\fontsize{16}{16}$\varepsilon$}
\psfrag{-e}{\fontsize{16}{16}$-\varepsilon$}
\psfrag{c}{\fontsize{16}{16}$v + \xi / 2$}
\psfrag{d}{\fontsize{16}{16}$v - \xi / 2$}
\psfrag{t1}{\fontsize{16}{16}$t_1$}
\psfrag{t2}{\fontsize{16}{16}$t_2$}
\psfrag{O}{\fontsize{16}{16}$O$}
\psfrag{A}{\fontsize{14}{14}$\om > \f 4 {\gamma^2} \f {\mu + 1} \mu$}
\psfrag{B}{\fontsize{14}{14}$\om = \f 4 {\gamma^2} \f {\mu + 1} \mu$}
\psfrag{C}{\fontsize{14}{14}$\om < \f 4 {\gamma^2} \f {\mu + 1} \mu$}
\psfrag{D}{\fontsize{14}{14}$\om_1$}
\psfrag{E}{\fontsize{14}{14}$\om_2$}

\newtheorem{theorem}{Theorem}[section]

\theoremstyle{definition}

\newtheorem{prop}[theorem]{Proposition}

\theoremstyle{remark}

\newcommand{\Hmu}{H_\mu^1}
\newcommand\vv{\textsc{v}}
\newcommand{\bes}{{\begin{split}}}
\newcommand{\ees}{{\end{split}}}
\newcommand{\bees}{{\begin{equation}\begin{split}}}
\newcommand{\es}{{\end{split}\end{equation}}}
\newcommand{\erre}{{\mathbb R}}
\newcommand{\grad}{\nabla}
\newcommand{\R}{{\mathbb R}}

\newcommand{\G}{{\mathcal G}}
\newcommand{\EE}{{\mathcal E}}
\newcommand{\C}{{\mathbb C}}

\newcommand{\bea}{\begin{eqnarray}}

\newcommand{\eea}{\end{eqnarray}}

\newcommand{\be}{\begin{equation}}

\newcommand{\ee}{\end{equation}}

\newcommand{\f}{\frac}

\newcommand{\ve}{\varepsilon}

\newcommand{\om}{\omega}

\newcommand{\ov}{\overline}

\numberwithin{equation}{section}

\newcommand{\blue}[1]{{\color{blue}#1}}
\newcommand{\red}[1]{{\color{red}#1}}

\tikzstyle{nodo}=[circle,draw,fill,inner sep=0pt,minimum size=\widthof{k}]
\tikzstyle{infinito}=[circle,inner sep=0pt,minimum size=0mm]
\tikzstyle{nodino}=[circle,draw,fill,inner sep=0pt,minimum size=0.5mm]

\begin{document}
\large

\title{Ground states for NLS on graphs: \\
a subtle interplay of metric and topology}


\author{Riccardo Adami}

\address{Dipartimento di Scienze Matematiche G.L. Lagrange,
Politecnico di Torino, Italy}




\date{today}


\keywords{}

\begin{abstract}
We review some recent results on the minimization of the energy
associated to the nonlinear Schr\"odinger Equation on non-compact
graphs.
Starting from seminal results given by the author together
with C. Cacciapuoti, D. Finco, and D. Noja for the star graphs, we
illustrate the achiements attained for general graphs and the related
methods, developed in collaboration with E. Serra and P. Tilli. 
We emphasize ideas and examples rather than computations or proofs.
\end{abstract}

\maketitle

\section{Introduction}
The subject of nonlinear evolution on {\em graphs} or {\em networks},
first introduced by F. Ali Mehmeti (\cite{ali}) and developed by
several authors (e.g. \cite{vonbelow,matrasulov,hannes,bellazzini,caudrelier}),   has rapidly become highly
popular in a quite spread scientific community, ranging from experts
in pointwise potentials (\cite{berkolaiko,exner,kuchment}) up to specialists of the
Nonlinear Schr\"odinger Equation and its standing waves (\cite{acfn-epl,cacciapuoti,kevrekidis,marzuola,NPS}).

The interest was initially driven by physical applications (for an
exhaustive introduction see \cite{noja14}), that involve
propagation in optical fibers and junctions  (\cite{bulgakov,smi,malomed}), 
and Bose-Einstein condensation (\cite{tokuno,vidal,zapata}). Subsequently, 
interesting mathematical issues also arose and have been considered as relevant for themselves.

Here we restrict our scope to the most basic variational problem: proving the
existence or the nonexistence of the ground
state for the Schr{\"o}dinger equation endowed with a focusing
cubic nonlinearity, on a
non-compact graph $\G$. The dynamics we investigate is then defined by the
equation
\be \label{schrod}
i \partial_t u = - \Delta u - |u|^{2} u,
\ee
but, since we concentrate on the problem of the ground state, the
central mathematical object in our analysis is the {\em energy functional}
(\cite{acfn-rmp,an-jpa})
\be \label{energy}
E (u) = \frac 1 2 \| \grad u \|^2_{L^2 (\G)} - \frac 1 4 \| u
\|^4_{L^4 (\G)},
\ee
whose value, as widely known, is conserved by the dynamics
generated by equation \eqref{schrod}.

Let us take for a moment the exact meaning of the functional
spaces involved in formula \eqref{energy} and of the Laplacian
appearing in \eqref{schrod} apart, and focus preliminarily on the physical
interpretation of the problem.
In fact, finding the ground state of the energy functional \eqref{energy} means
finding the wave function of a Bose-Einstein condensate made of atoms
attracting one another, placed in an
optical and/or magnetic trap supplied with many branches and ramifications. The
non-compactness of the graph is the mathematical translation of the fact that some of the
branches of the trap are much longer than others, and also much longer
than a characteristic length of the condensate, to be specified later.

We remark that the expression of the energy \eqref{energy} is formally equivalent to the 
analogous formula
for the NLS in $\R^N$. The non standard elements here are the graph $\G$, and,
therefore, the functions defined on it and the related functional
spaces.
However, the
definition of $\G$ as a {\em metric} graph, of the functions defined on
it, and of the associated $L^q$-spaces, are the most natural ones. We rapidly
review them for the convenience of the reader.

A {\em graph} $\G$ is defined by two sets: a set $V$ of {\em
  vertices}, that are to be thought of as spatial points, and a set $E$ of
{\em edges} (or {\em bonds}), each of them to be thought of as a line
joining a couple of
vertices. It is convenient to isolate a subset of $V$, denoted by
$V_\infty$, that contains the {\em vertices at infinity}. We establish
that two distinct vertices
at infinity cannot be connected by an edge. Furthermore, defining the
{\em degree} of a vertex as the number of edges
starting or ending on it, we assume that every vertex at infinity
has degree one, i.e. it is the endpoint of exactly one
halfline. We also notice that there is no loss of generality in
supposing that every vertex has degree at least three, and this is
what we shall always do.

\noindent
For the sake of simplicity, we require the resulting structure
to be {\em connected}, i.e. starting from an arbitrary vertex it is
possible to reach any other vertex by following a sequence of adjacent edges.

By all this,
the {\em topology} of the graph is fixed. Let us turn to the {\em metric}.

\noindent
In order to endow $\G$ with a metric structure, we identify
each edge $e$ with a real interval $I_e : = [0, \ell_e]$, with $\ell_e > 0$, or,
if the edge ends up into a vertex at infinity,
with the halfline $I_e := [0, + \infty)$. Thus, $\G$ is a metric space formed
  by the union of the intervals representing the edges, and the
  distance between two points is given by the shortest path joining
  the two points through a
  connected sequence of edges.
In this way we obtain a {\em metric graph} (Fig.1).
\begin{figure}[t]
\begin{center}
\begin{tikzpicture}[xscale= 1.5,yscale=1.5]
%
\node at (0,0) [nodo] (1) {};
\node at (-1.5,0) [infinito]  (2){$\infty$};
\node at (1,0) [nodo] (3) {};
\node at (0,1) [nodo] (4) {};
\node at (-1.5,1) [infinito] (5) {$\infty$};
\node at (2,0) [nodo] (6) {};
\node at (3,0) [nodo] (7) {};
\node at (2,1) [nodo] (8) {};
\node at (3,1) [nodo] (9) {};
\node at (4.5,0) [infinito] (10) {$\infty$};
\node at (5.5,0) [infinito] (11) {$\infty$};
\node at (4.5,1) [infinito] (12) {$\infty$};
\draw [-] (1) -- (2) ;
 \draw [-] (1) -- (3);
 \draw [-] (1) -- (4);
 \draw [-] (3) -- (4);
 \draw [-] (5) -- (4);
 \draw [-] (3) -- (6);
 \draw [-] (6) -- (7);
 \draw [-] (6) to [out=-40,in=-140] (7);
\draw [-] (3) to [out=10,in=-35] (1.4,0.7); 
\draw [-] (1.4,0.7) to [out=145,in=100] (3); 
 \draw [-] (6) to [out=40,in=140] (7);
 \draw [-] (6) -- (8);
 \draw [-] (6) to [out=130,in=-130] (8);
 \draw [-] (7) -- (8);
 \draw [-] (8) -- (9);
  \draw [-] (7) -- (9);
  \draw [-] (9) -- (12);
  \draw [-] (7) -- (10);
  \draw [-] (7) to [out=40,in=140] (11);
\end{tikzpicture}
\end{center}
\caption{\footnotesize{{\bf A typical graph.} Vertices are bullets, edges are lines, vertices
at infinity are denoted by $\infty$. Notice that self-loops and multiple connections are allowed.}}
\end{figure}
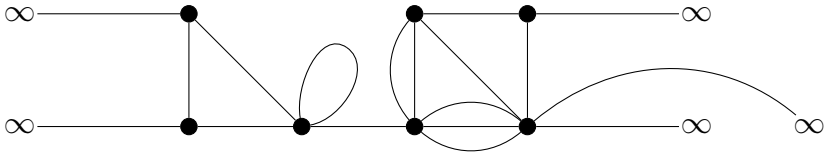

As a next step, the notion of function (or, thinking of quantum
mechanical applications, {\em wave function}) on the metric space $\G$
is natural: the function $u : \G \longrightarrow \C$ is defined as the set of
functions $u_e, \, e \in E$, i.e. $u_e$ is the restriction of $u$ to the
edge $e$. For a function $u$ to be continuous, in addition to the
continuity of every restriction
$u_e$, one has to require {\em continuity at vertices}: for instance, if the same
vertex $\vv$ is the initial point of two edges $e$ and $e'$,
then  $u_e (0) = u_{e'} (0)$, and an analogous equality holds if the
same vertex acts as initial point for a vertex $e$ and as endpoint for
$e'$. In the same way we define differentiability, and denote by
$\grad u$ the complex-valued function on $\G$ whose restriction to the
edge $e$ is given by $u'_e$.

\noindent
Once we have introduced functions on $\G$, it is immediate to define
functional spaces. We need  the spaces $L^q
(\G)$, defined by 
$$L^q (\G) = \oplus_{e \in E} L^q (I_e), \qquad \| u
\|^q_{L^q (\G)} \ = \ \sum_{e \in E} \| u_e \|^q_{L^q (I_e)},$$ and the
space $H^1 (\G)$, defined by   
$$H^1 (\G) = \oplus_{e \in E} H^1 (I_e), \qquad \| u
\|^2_{H^1 (\G)} \ = \ \sum_{e \in E} \| u_e \|^2_{H^1 (I_e)}$$
 with the
additional requirement of the continuity of $u$ (i.e. $u$ must be
continuous at every vertex).

Thus, the energy functional \eqref{energy} can be rewritten as
$$
E (u, \G) \ = \ \frac 1 2 \sum_{e \in E} \int_0^{\ell_e} |u'_e|^2 dx -
\frac 1 p  \sum_{e \in E} \int_0^{\ell_e} |u_e|^p dx,
$$
where, with a slight abuse of notation, we allowed $\ell_e = + \infty$,
and as a domain for $E (\cdot, \G)$ we choose $H^1 (\G)$. Notice that
this is not the only possible option: another suitable domain could be $\oplus_{e
  \in E} H^1 (I_e)$. In that case, one would find that minimizers
  satisfy homogeneous Neumann
conditions at every vertex instead of continuity: obviously, the
ground state would change, and in general both existence and shape of the ground states, as well
as of the stationary waves, strongly depend on conditions at the vertices, and therefore
on the choice of the enery domain.

It is immediately seen that there is no absolute minimum
of $E(u, \G)$:
indeed, for any non-vanishing $u$, one has
$E (\lambda u, \G) \to - \infty$ as $\lambda \to + \infty$, i.e. for
large data the energy is unbounded from below. However, since
the flow associated to the NLS preserves the $L^2$-norm, or {\em
  mass}, it is physically meaningful to seek minimizers of the energy functional
\eqref{energy} under the constraint of {\em constant mass}. Indeed,
in current experiments on Bose-Einstein condensates, what we are
calling the mass denotes
the number of atoms (in different contexts related to quantum physics,
it denotes
the overall probability of finding a
particle on a certain region and so on), therefore the attainable
configuration of minimal energy is always conditioned by the choice of
the initial mass and this identifies a natural constraint.

Then, the problem we are interested in is the following.

\medskip

\noindent
{\bf Problem P.} {\em Given $\mu > 0$, defined the space
$$
H^1_\mu (\G) = \{ u \in H^1 (\G), \, \| u \|_{L^2(\G)} = \mu \},
$$
and introduced the notation
$$
\mathcal E_\G (\mu) \ : = \ \inf_{u \in H^1_\mu (\G)} E (u, \G),
$$
find a function $u \in H^1_\mu (\G)$ such that
$$E (u, \G) = \mathcal E_\G (\mu).$$}

\medskip

Notice that this is a standard variational problem (on a non-standard
environment),
and involves two
main elements: first, the {\em boundedness from below} of the {\em
  constrained} energy functional; second, the existence of a function
  $u$ whose energy equals the infimum of the constrained
functional. Such a function is also called a {\em minimizer}.

As a preliminary fact we stress that, despite the non-boundedness of
the energy functional \eqref{energy} on the domain $H^1 (\G)$, the
constrained energy is bounded from below:
indeed, Gagliardo-Nirenberg's type estimates hold for
graphs too \cite{ast-arxiv,t16}, in particular
\be \label{gagliardo}
\| u \|_{L^p (\G)}^p \ \leq \ C \| u \|_{L^2 (\G)}^{\frac p 2 + 1}  \|
\nabla u \|_{L^2 (\G)}^{\frac p 2 - 1},
\ee
that, applied to the energy functional \eqref{energy}, gives
\be \label{lowerbound}
E (u, \G) \ \geq \ \f 1 2 \| \nabla u \|_{L^2 (\G)}^2 \left( 1 - C
\frac {\mu^{\frac 3 2}}{\|\nabla u  \|_{L^2 (\G)}}
\right)
\ee
which, for fixed $\mu$, is a lower bounded quantity, so that
\be \label{lb}
\inf_{u \in H^1_\mu} E (u, \G) \ > \ - \infty.
\ee
It remains then to establish whether the infimum is attained or not.
In order to suitably review the result on that, it is convenient to
recall some preliminary notions.

\subsection{Preliminary results}

Since $E (|u|, \G) \leq E (u, \G)$, from now on we shall restrict to
real non-negative functions.

First, a ground state $u$ gives a {\em standing wave} for the
dynamical problem \eqref{schrod}, in the sense that the function
$$
\psi (t,x) \ = \ e^{i \om t} u (x)
$$
is a solution to \eqref{schrod}. Of course, plugging the expression of
$\psi$ into equation \eqref{schrod} gives
\be \label{lagrange}
\Delta u + u^3 \ = \ \omega u,
\ee
where $\Delta u$ is the function whose restriction on the edge $e$ is
$u_e''$, and the requirement that $u$ is a stationary point for the
functional $E (\cdot, \G)$ provides the so-called Kirchhoff condition
at vertices, namely
\be \label{kirchhoff}
\sum_{e\succ \vv}
\frac{d u_e}{d x_e}(\vv)=0
\ee
where $e\succ \vv$ denotes that the edge $e$ is incident at the vertex
$\vv$. Of
course, there can be stationary solutions to \eqref{schrod} that  are
not ground states, and they are usually called {\em excited states.}

\medskip

It is well-known that the version on the line of Problem P, i.e. when $\G =
\R$, is solved by the {\em soliton} (see Fig.2)
$$
\phi_\mu (x) \ = \ \f \mu {2 \sqrt 2} \, \hbox{{sech}} \left( \f \mu 4
x \right).
$$
and by its translated. Notice that the existence and the uniqueness of
the soliton at mass $\mu$ fixes a natural lengthscale on the graph as
$\mu^{-1}$. There is then a correpondence between {\em large mass} and
{\em long edges}, rooted in the fact that a soliton with a large mass
is squeezed, and thus it sees edges as long.
\begin{figure} \label{soliton}
\begin{center}
\begin{tikzpicture}[xscale= 0.9,yscale=0.9]
\node at (0,0) [nodo]{};
\draw[-,thick,blue] (0,2) to [out=0,in=170] (2.5,0.2);
\draw[-,thick,blue] (0,2) to [out=180,in=10] (-2.5,0.2);
\draw[-,thick,blue] (2.5,0.2) to [out=350,in=178] (6,0);
\draw[-,thick,blue](-2.5,0.2) to [out=190,in=2] (-6,0);
\draw[-] (-6,0) -- (6,0);
\node at (2.2,1.5) {$\blue{\phi_\mu}$};
\end{tikzpicture}
\end{center}
{\caption {\footnotesize{{\bf The soliton $\phi_\mu$}. It is the
      unique ground state on the line, up to translations and
      multiplications by a phase factor.}}}
\end{figure}
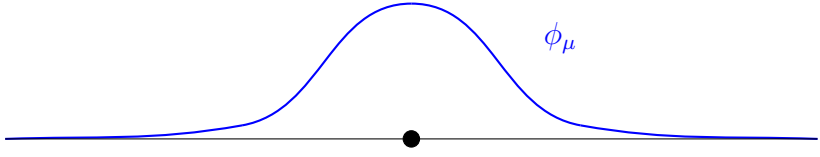
The energy level of $\phi_\mu$ reads
$$
\EE_\R(\mu) \ = - \f {\mu^{3}}{96}
$$
and plays a crucial role in Problem P. Indeed, being the graph $\G$
non-compact, it contains at least one halfline, where it is always
possible to approximate arbitrarily well a soliton. So that,
\be \label{first}
\EE_\G (\mu) \ \leq  \ - \f {\mu^{3}}{96}.
\ee
On the other hand, it is clear that for some graphs it is possible to
go below the level $\EE_\R (\mu)$. For instance, if $\G = \R^+$,
then the ground states at mass $\mu$ are obviously given by half
solitons $\chi_+ \phi_{2 \mu}$
centred at the origin
of the halfline, and therefore
$$
\EE_{\R^+} (\mu) \ = \ - \f {\mu^3} {24}.
$$
The level $\EE_{\R^+} (\mu)$ of the halfline plays an important role
too,
as it represents the minimal possible level among all non-compact
graphs. In order to understand this point, one must introduce the
following result:
\begin{prop}[Monotone rearrangement] \label{m-rearrangement}
Given $u \in H^1_\mu (\G)$, there exists a function $u^* \in H^1_\mu
(\R^+)$ s.t.
$$ \| u^* \|_{L^4 (\R^+)} =  \| u \|_{L^4 (\G)} , \qquad
\| (u^*)' \|_{L^2 (\R^+)} =  \| \nabla u \|_{L^2 (\G)}
$$
\end{prop}
In practice, the function $u^*$ can be constructed as the {\em
  monotone rearrangement} of $u$ on the halfline (see
\cite{ast-cv,friedlander}),
and it is clear that
\be \label{second}
E (u^*, \R^+) \ \leq \ E (u, \G).
\ee
The results in \eqref{first} and \eqref{second} can be put together to
state the following proposition:
\begin{prop}[Pinching] \label{pinching}
For every non-compact graph $\G$, the following estimates
hold:
$$
- \f {\mu^3} {24} \ \leq \ \EE_\G (\mu) \ \leq \ - \f {\mu^3} {96}
$$
\end{prop}
Another important preliminary result is the following:
\begin{prop} [Comparison] \label{comparison}
The infimum of $E (\cdot, \G)$ in $H^1_\mu (\G)$ is attained if and
only if
there exists $u \in H^1_\mu (\G)$ such that
$$ E (u, \G) \ \leq \ - \f {\mu^3} {96}. $$
\end{prop}
Such a proposition is relevant under both  a practical and a conceptual point of view. At a practical level, it shows that in order to prove existence of the ground state, it is sufficient to exhibit a function whose energy level is lower than the level of the soliton. Conceptually, first,
as already explained, there are functions whose energy is arbitrarily
close to that of the soliton, but, except for some particular graphs,
a soliton cannot be exactly constructed, so that the energy level of a
soliton may in general not be attained, because sequences that
approximate a soliton run away at infinity, weakly converging to zero.

\noindent
Now, Proposition \ref{comparison} establishes that if there is a
function whose energy is lower
than the energy of the soliton with the same mass, then there are
minimizing sequences that are compact, i.e., they do not run away in
order to reconstruct a far away soliton, but converge to a minimizer.

In other words, there is a competition between the soliton and the
function that can be hosted on a graph. If the latter win the
competition, then the infimum is attained.

\medskip
In what follows we shall give examples and results that illustrate how topology
(Section \ref{sec:topology}) or metric (Section \ref{sec:metric})
may influence the existence of a ground state.


\section{Topology} \label{sec:topology}

As mentioned in the previous section, the halfline can be seen as a
graph with one vertex at the origin and the other at infinity. On the
other hand, the line can be understood as a graph made of two
halflines and three vertices: one joining the two halflines, the two
others at infinity.
These two cases can also be described as {\em infinite star graphs},
i.e. graphs made of a certain number of halflines, meeting one another
at a unique vertex.  Immediately beyond the case of the halfline and
of the line, lies the case of the star graph $\mathcal S_3$ made of {\em three}
halflines connected at a single vertex (Fig.3):
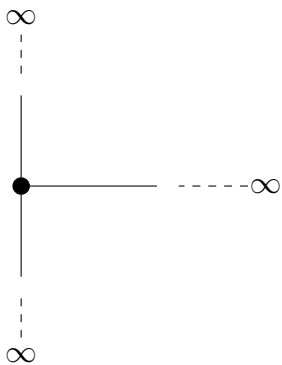
\begin{figure}[h]
\begin{center}
\begin{tikzpicture}[xscale= 1.3,yscale=1.5]
\node at (0,0) [nodo] (2) {};
\node at (2.5,0) [infinito]  (3) {$\infty$};
\node at (0,1.5) [infinito]  (4) {$\infty$};
\node at (0,-1.5) [infinito]  (8) {$\infty$};
\node at (1.5,0)  [minimum size=0pt] (5) {};
\node at (-1.5,0) [minimum size=0pt] (6) {};
\node at (0,0.9) [minimum size=0pt] (7) {};
\node at (0,-0.9) [minimum size=0pt] (9) {};

\draw[dashed] (3) -- (5);
\draw[dashed] (7) -- (4);
\draw[dashed] (9) -- (8);
\draw [-] (5) -- (2) ;
\draw [-] (2) -- (7) ;
\draw [-] (2) -- (9) ;
\end{tikzpicture}
\end{center}
{\caption {\footnotesize{{\bf The infinite three-star graph $\mathcal S_3$}.
Three halflines meeting one another at a single vertex. For this
graph, whatever the value of the mass $\mu$, there is no ground state, as
minimizing sequences run away towards infinity mimicking the shape of
the soliton $\phi_\mu$ on a single halfline.}}}
\end{figure}

\noindent
It has been shown in \cite{acfn-jpa} (see also \cite{acfn-aihp,acfn-jde}) that in this case
$$
\EE_{\G} (\mu) \ = \ - \f{\mu^3}{96}
$$
but such an infimum is not attained. The situation can be displayed as
follows: the
only stationary state $u_s$ is made of three half-solitons of mass $\mu / 3$
with their maximum at the vertex (Fig.4). 
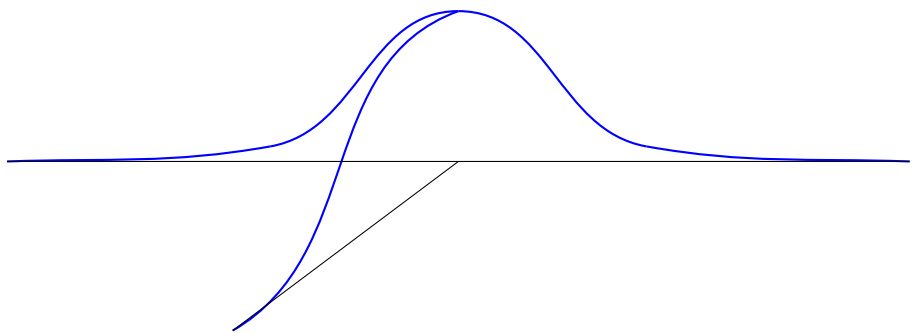
\begin{figure} \label{three}
\begin{center}
\begin{tikzpicture}[xscale= 1,yscale=1]
\draw[-,thick,blue] (0,2) to [out=0,in=170] (2.5,0.2);
\draw[-,thick,blue] (0,2) to [out=180,in=10] (-2.5,0.2);
\draw[-,thick,blue] (2.5,0.2) to [out=350,in=178] (6,0);
\draw[-,thick,blue] (-2.5,0.2) to [out=190,in=2] (-6,0);
\draw[-,thick,blue] (0,2) to [out=200,in=30] (-3,-2.25);
\draw[-] (-6,0) -- (6,0);
\draw[-] (-3,-2.25) -- (0,0);
\end{tikzpicture}
\end{center}
\caption{\footnotesize{{\bf Three half-solitons.} This is the only stationary state for the
infinite star-graph made of three half-line. However, it is not a ground state.}}
\end{figure}

The energy level of $u_s$ can be
easily computed as
$$
E (u_s, \G) \ = \ - \frac {\mu^3} {216} \ > \ - \frac {\mu^3} {96}
$$
so that $u_s$ cannot be a ground state, in spite of the fact that it
is the only stationary state, i.e. the only solution to equation
\eqref{lagrange}.
The same phenomenon occurs for all star-graphs made of $N \geq 3$
halflines, as outlined in \cite{acfn-jpa}. The deep reason for that has
been later investigated in \cite{ast-cv}, and turns out to be rooted in the
rearrangement theory.

\subsection{Rearrangements: counting preimages}
The monotone rearrangement $u^*$, introduced in Proposition \ref{m-rearrangement},
is not the unique rearrangement that can be made in order to pass from
functions on the graph to functions on the line. Another kind of
rearrangement is the {\em symmetric} one, that associates to any $u
\in H^1_\mu (\G)$ a function $\widehat u \in  H^1_\mu (\R)$. The following
theorem holds:
\begin{prop}[Symmetric rearrangement] \label{ms-rearrangement}
{\em Given $u \in H^1_\mu (\G)$, there exists a function $\widehat u \in H^1_\mu
(\R)$ s.t.
$$ \| \widehat u \|_{L^4 (\R)} =  \| u \|_{L^4 (\G)}.$$
Furthermore, if almost every element of Ran $u$ has at least two
preimages in the domain of $u$, then
$$
\| (\widehat u)' \|_{L^2 (\R)} \leq  \| \nabla u \|_{L^2 (\G)}.
$$}
\end{prop}
For practical purposes, $\widehat u$ can be constructed as the unique
function on $\R$, that shares with $u$ the same distribution function
and is
symmetric and monotonically decreasing in $[0, + \infty)$. The
  relationship with the monotone rearrangement writes simply as
$$ \widehat u (x) : = u^* (2 |x|).$$

Roughly speaking, Proposition \ref{ms-rearrangement} shows that
the rearranged function makes the kinetic energy decrease as it lowers
the  measure of the set of the preimages of almost all elements of Ran
$u$:
being the set of preimages a zero-dimensional manifold, the
natural measure is given by counting points. In this respect, the
monotonic rearrangement $u^*$ introduced in the previous section is the best
one could find, since for any element in Ran $u$ it provides {\em one}
preimage only. On the other hand, the symmetric rearrangement
$\widehat u$ provides {\em two} preimages for every element of Ran
$u$ (except for the maximum of $u$), so that symmetrically rearranging
lowers the energy under the additional hypothesis that for almost
every point in Ran $u$ there are at least two preimages.

We are now in position to understand the reason of the negative result
for the infinite star graph made of three halflines: consider
a function $u \in H^1_\mu (\mathcal S_3)$, and imagine to follow
the function starting from its maximum point: you can move rightwards
or leftwards (unless the maximum is not at the vertex, where one has
possibly still more options), but whatever direction you choose, in
order to belong to $H^1 (\mathcal S_3)$ the value of the function you are
following must end up at
zero. In other words, for every point in Ran $u$ (except at most for
max $u$),
there are at least
two preimages, one for each direction that can be taken starting from
the maximum of $u$, so that
$$
E (u, \mathcal S_3) \ \geq \ E (\widehat u, \R) \ \geq \ E (\phi_\mu, \R) \ =
\ - \f {\mu^3} {96}.
$$
Furthermore, in order to attain the level $ - \f {\mu^3} {96}$, the
rearrangement $\widehat u$ must coincide with $\phi_\mu$, then, owing to the construction of
the rearrangement $\widehat u$ (see \cite{ast-cv}), the support
of $u$ should be infinite, so that $u$ can be considered as a function on the line. As a consequence, in order to avoid having
more than two preimages in an interval of elements of its range, $u$
should be entirely supported in a halfine and attain zero at the
origin, so $u$
 would have a {\em
  strictly greater energy} than the soliton, contradicting the
assumption that the level of the soliton is attained. This is why
there is no ground state if $\G$ is an infinite star-graph with more
than two halflines, regardless of the value of the mass.

\subsection{Assumption (H): graphs as bunches of lines}
It is worth remarking that, for an infinite star graph
with more than two halflines,
the impossibility of having less than two preimages is dictated by the
topological structure of the graph: metrics has no role when one
considers all functions in $H^1_\mu (\G)$, however shrunk or
stretched. More specifically, in our argument the crucial fact was
that, starting from the maximum point of $\G$, a point at infinity was
available in both directions, so that, in order to belong to $H^1
(\G)$, the function $u$ was forced to go to zero in both directions,
running through its values at least twice. The same necessity is
not present in the case of the halflines, where, if one
places the maximum at the origin, then it is possible to attain each
value in the range of $u$ exactly once.

In the effort of understanding the relationship between
topology of a graph and existence of ground states, a non-trivial task then was involved in the search for topological
conditions that ensure that

\medskip

\noindent
{\em Given an arbitrary function $u \in H^1
  (\G)$, almost every point in Ran $u$ is endowed with at least two
  preimages}.

\medskip

The key point is to find a simple condition capable to ensure that, when following a function starting from its
maximum, a point at infinity is available in both
directions. Besides, since the maximum may be located at any point
of the graph, one has to guarantee that, starting from {\em every} point on the
graph, an infinity point is available when moving in either
direction.

We finally got the following condition (Fig.5):

\medskip
\noindent
{\bf Assumption (H). First version}. {\em Every point of $\G$ lies in
a trail that
  connects two different vertices at infinity.}

\medskip
A trail is a connected sequence of non-repeated edges. Assumption (H)
states exactly that from every point one can get to infinity through
two disjoint paths, so the argument used for the infinite star graphs
still holds. In other words, every edge can be considered as part of a
line, so it is not possible to get an energy level lower than the energy of
the soliton.
\begin{figure}[t]
\begin{center}
\begin{tikzpicture}[xscale= 1.4,yscale=1.4]
\node at (0,0) [nodo] (1) {};
\node at (-1.5,0) [infinito]  (2){$\infty$};
\node at (1,0) [nodo] (3) {};
\node at (0,1) [nodo] (4) {};
\node at (-1.5,1) [infinito] (5) {$\infty$};
\node at (2,0) [nodo] (6) {};
\node at (3,0) [nodo] (7) {};
\node at (2,1) [nodo] (8) {};
\node at (3,1) [nodo] (9) {};
\node at (4.5,0) [infinito] (10) {$\infty$};
\node at (5.5,0) [infinito] (11) {$\infty$};
\node at (4.5,1) [infinito] (12) {$\infty$};
\node at (1.5,0.85) {$\scriptstyle x$};
\node at (1.4,0.7) [nodino]  {};

\draw [-] (1) -- (2) ;
\draw [-] (1) -- (3);
\draw [-] (1) -- (4);
\draw [-] (3) -- (4);
\draw [-] (5) -- (4);
\draw [] (3) -- (6);
\draw [-] (6) -- (7);
\draw [-] (6) to [out=-40,in=-140] (7);
\draw [-] (3) to [out=10,in=-35] (1.4,0.7);
\draw [-] (1.4,0.7) to [out=145,in=100] (3);
\draw [-] (6) to [out=40,in=140] (7);
\draw [-] (6) -- (8);
\draw [-] (6) to [out=130,in=-130] (8);
\draw [-] (7) -- (8);
\draw [-] (8) -- (9);
\draw [-] (7) -- (9);
\draw [-] (9) -- (12);
\draw [-] (7) -- (10);
\draw [-] (7) to [out=40,in=140] (11);
\end{tikzpicture}
\end{center}
\end{figure}

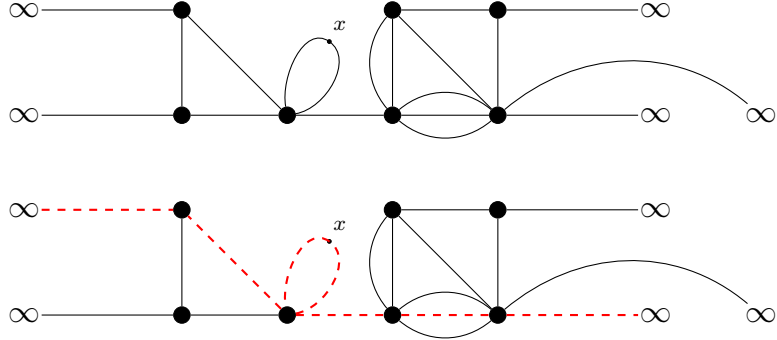
\begin{figure}[t]
\begin{center}
\begin{tikzpicture}[xscale= 1.4,yscale=1.4]
\node at (0,0) [nodo] (1) {};
\node at (-1.5,0) [infinito]  (2){$\infty$};
\node at (1,0) [nodo] (3) {};
\node at (0,1) [nodo] (4) {};
\node at (-1.5,1) [infinito] (5) {$\infty$};
\node at (2,0) [nodo] (6) {};
\node at (3,0) [nodo] (7) {};
\node at (2,1) [nodo] (8) {};
\node at (3,1) [nodo] (9) {};
\node at (4.5,0) [infinito] (10) {$\infty$};
\node at (5.5,0) [infinito] (11) {$\infty$};
\node at (4.5,1) [infinito] (12) {$\infty$};
\node at (1.5,0.85) {$\scriptstyle x$};
\node at (1.4,0.7) [nodino]  {};

\draw [-] (1) -- (2) ;
\draw [-] (1) -- (3);
\draw [-] (1) -- (4);
\draw [thick,red,dashed] (3) -- (4);
\draw [thick,red,dashed] (5) -- (4);
\draw [thick,red,dashed] (3) -- (6);
\draw [thick,red,dashed] (6) -- (7);
\draw [-] (6) to [out=-40,in=-140] (7);
\draw [thick,red,dashed] (3) to [out=10,in=-35] (1.4,0.7);
\draw [thick,red,dashed] (1.4,0.7) to [out=145,in=100] (3);
\draw [-] (6) to [out=40,in=140] (7);
\draw [-] (6) -- (8);
\draw [-] (6) to [out=130,in=-130] (8);
\draw [-] (7) -- (8);
\draw [-] (8) -- (9);
\draw [-] (7) -- (9);
\draw [-] (9) -- (12);
\draw [thick,red,dashed] (7) -- (10);
\draw [-] (7) to [out=40,in=140] (11);
\end{tikzpicture}
\end{center}
\caption{\footnotesize{{\bf Assumption (H):} every point of the graph lies on a
  trail connecting two distinct vertices at infinity. The point $x$
  belongs to several such trails: in the picture we put in evidence
  one of them, by representing it through dashed lines.}}
\end{figure}









%




An equivalent version of Assumption (H), that proves more manageable
for proofs, is the following:

\medskip

\noindent
{\bf Assumption (H). Second version}. {\em Removing an arbitrary edge,
  every resulting connected component contains a vertex at infinity.}

\medskip
\noindent
Owing to this version, it is easier to visualize how Assumption (H)
can fail. Clearly, it fails when there is one halfline only (for
instance, if $\G = \R^+$), and it fails when there is a terminal
edge, or ``pendant'' too:
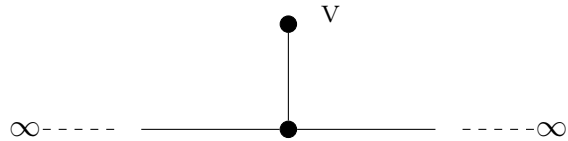
\begin{figure}[h]
\begin{center}
\begin{tikzpicture}[xscale= 1.4,yscale=1.4]
\node at (-2.5,0) [infinito]  (1) {$\infty$};
\node at (0,0) [nodo] (2) {};
\node at (2.5,0) [infinito]  (3) {$\infty$};
\node at (1.5,0)  [minimum size=0pt] (5) {};
\node at (-1.5,0) [minimum size=0pt] (6) {};
\node at (0,1) [nodo] (7) {};
\node at (0.4,1.1) {$\vv$};

\draw[dashed] (1) -- (6);
\draw[dashed] (3) -- (5);
\draw [-] (2) -- (6) ;
\draw [-] (5) -- (2) ;
\draw [-] (2) -- (7) ;
\end{tikzpicture}
\end{center}
\end{figure}

\begin{figure}[h] \label{2+0}
\begin{center}
\begin{tikzpicture}[xscale= 1.4,yscale=1.4]
\node at (-2.5,0) [infinito]  (1) {$\infty$};
\node at (0,0) [nodo] (2) {};
\node at (2.5,0) [infinito]  (3) {$\infty$};
\node at (1.5,0)  [minimum size=0pt] (5) {};
\node at (-1.5,0) [minimum size=0pt] (6) {};
\node at (0,1) [nodo] (7) {};
\node at (0.4,1.1) {$\vv$};

\draw[dashed] (1) -- (6);
\draw[dashed] (3) -- (5);
\draw [-] (2) -- (6) ;
\draw [-] (5) -- (2) ;
\end{tikzpicture}
\end{center}
{\caption {\footnotesize{{\bf Failure of Assumption (H) for the line
with a pendant}.
Removing the pendant, two connected components are produced: one of
them is a line, the second one is made by a vertex at a finite point
of the space. So the line with a pendant, although endowed with two
halflines, violates Assumption (H).}}}
\end{figure}
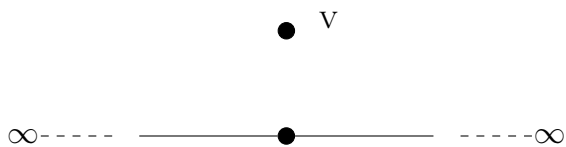
Notice that, in the latter case, violation of Assumption (H) is due to the
presence of the pendant, whose removal yields a compact connected component (Fig.6).

\medskip
Assumption (H) leads to a negative result, that is the main goal of
\cite{ast-cv}.
\begin{theorem}[Nonexistence] \label{nonexistence}
Assume that $\blue\G$ satisfies assumption \blue{(H)}. Then
$$
\inf_{u \in H^1_\mu ({\mathcal G})} E (u, {\mathcal G}) \ = \ E
(\phi_\mu, \erre)
$$
and it is {never attained}, except if $\G$
is a ``bubble tower''.

\end{theorem}

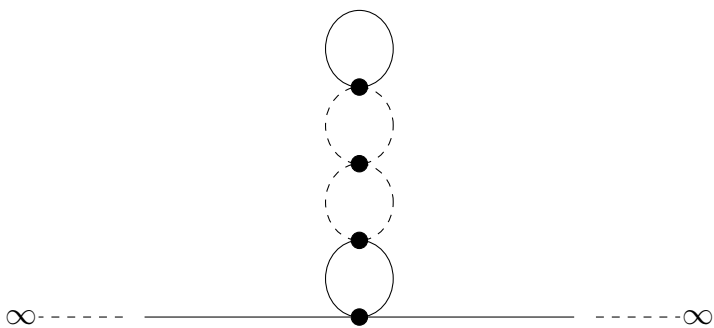
\begin{figure}[h] \label{e3}
\begin{center}
\begin{tikzpicture} [xscale= 1.5,yscale=1.7]
\node at (-3,0) [infinito]  (1) {$\infty$};
\node at (0,0) [nodo] (2) {};
\node at (3,0) [infinito]  (3) {$\infty$};
\node at (2,0)  [minimum size=0pt] (5) {};
\node at (-2,0) [minimum size=0pt] (4) {};
\node at (0,1.2) [nodo] (6) {};
\node at (0,1.8) [nodo] (7) {};
\node at (0,0.6) [nodo] (8) {};

\draw[dashed] (4) -- (1);
\draw [-] (4) -- (2);
\draw [-] (5) -- (2) ;
\draw[dashed] (5) -- (3);
\draw(0,0.3) circle (0.3cm);
\draw[dashed]  (0,0.9) circle (0.3cm);
\draw[dashed]  (0,1.5) circle (0.3cm);
\draw (0,2.1) circle (0.3cm);
\end{tikzpicture}
\end{center}
{\caption {\footnotesize{{\bf Bubble tower}.
Bubble towers are the only graphs that satisfy Assumption (H) but
admit ground states (again, regardless of the choice of $\mu$). Notice
that every bubble is cut in two arcs of equal length by the vertices
lying on it.}}}
\end{figure}
\medskip

The only graphs satisfying assumption (H) and admitting a minimum are
the ``bubble towers'', since one can cut a soliton on the line and
paste it on the tower (see Fig.7).
\begin{figure}[h]  \label{bubbles}  
\begin{center}
\begin{tikzpicture}[xscale= 0.7,yscale=1]
\draw[-,thick,blue] (-6,2) to [out=0,in=170] (-3.5,0.2);
\draw[-,thick,blue] (-6,2) to [out=180,in=10] (-8.5,0.2);
\draw[-,thick,blue] (-3.5,0.2) to [out=350,in=178] (-1,0);
\draw[-,thick,blue] (-8.5,0.2) to [out=190,in=2] (-11,0);
\draw[-] (-11,0) -- (-1,0);
\draw[dashed] (-8.4,0.5)--(-8.4,-0.5);
\draw[dashed] (-7.4,1.2)--(-7.4,-0.5);
\draw[dashed] (-4.6,1.2)--(-4.6,-0.5);
\draw[dashed] (-3.6,0.5)--(-3.6,-0.5);
\node at (-6,2.3)  {$u$};
\node at (-4,.8)  {$v$};
\node at (-8,.8)  {$v$};
\node at (-2.4,.4)  {$w$};
\node at (-9.6,.4)  {$w$};
\node at (-6,-0.3) {$\scriptstyle{\ell_1}$};
\node at (-4.1,-0.3) {$\scriptstyle{\ell_2/ 2}$};
\node at (-7.88,-0.3) {$\scriptstyle{\ell_2/ 2}$};
\draw[-] (0,0) -- (7,0);
\node at (3.5,0) [nodo] (1) {};
\node at (3.5,1.6) [nodo] (2) {};
\draw(3.5,.8) circle (.8cm);
\draw(3.5,2.2) circle (.6cm);
\draw [->] (2.5,1.1) arc [radius=.7, start angle=160, end angle= 200];
\draw [->] (4.5,1.1) arc [radius=.7, start angle=20, end angle= -20];
\draw [->] (2.7,2.5) arc [radius=.5, start angle=150, end angle= 210];
\draw [->] (4.3,2.5) arc [radius=.5, start angle=30, end angle= -30];
\draw[->] (2,0.2) -- (1.5,0.2);
\draw[->] (5,0.2) -- (5.5,0.2);
\node at (4,2.9) {$\scriptstyle{\ell_1}$};
\node at (4.4,1.5) {$\scriptstyle{\ell_2}$};
\node at (3,2.9) {$u$};
\node at (2.8,1.5) {$v$};
\node at (2,-.3) {$w$};
\node at (5,-.3) {$w$};
\end{tikzpicture}
{\caption {\footnotesize{{\bf A soliton on the line with two bubbles}.
A soliton can be
    cut and paste again on a bubble tower. Its energy is not affected
    by the procedure, so that the infimum is attained.}}}
\end{center}
\end{figure}
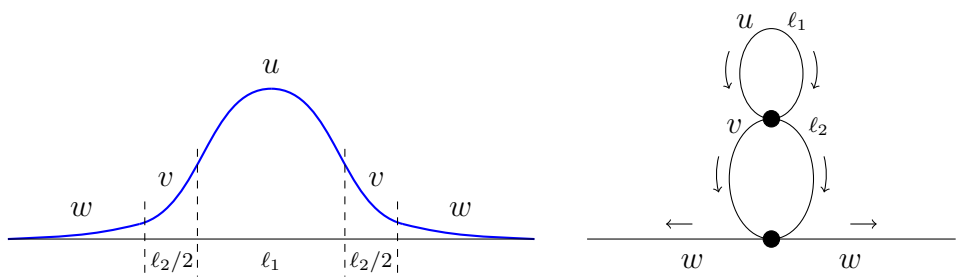
As the energy is not affected by the placement of the function on
the graph (Fig.8), Theorem
\ref{comparison} guarantees the existence of a ground state.

\subsection{Violating (H): a line with a pendant}
As already stated, one of the simplest
graphs that does not satisfy assumption (H) is the real line with a
pendant. We show that, for every $\mu > 0$, there exists a ground
state. The preocedure will also highlight the shape of such a ground
state. In particular, we prove that, denoted by $\mathcal P$ the graph
made of a line and a pendant,
$$
\inf_{u\in \Hmu(\G)} E(u,\mathcal P) < - \f  {\mu^3} {96},
$$
so that, owing to Proposition \ref{comparison}, the infimum is {attained}.
To this aim, We use a {\em graph surgery} together with rearrangements.

\medskip

\noindent
As a first step, we cut the soliton ${\phi_\mu}$ centred at a
  width ${\ell}$ (see Fig.9).

\begin{figure} \label{surg1}
\begin{center}
\begin{tikzpicture}[xscale= 1,yscale=1]
\draw[-,thick,blue] (0,2) to [out=0,in=170] (2.5,0.2);
\draw[-,thick,blue] (0,2) to [out=180,in=10] (-2.5,0.2);
\draw[-,thick,blue] (2.5,0.2) to [out=350,in=178] (6,0);
\draw[-,thick,blue] (-2.5,0.2) to [out=190,in=2] (-6,0);
\draw[-] (-6,0) -- (6,0);
\draw[red] (-4,1) -- (4,1);
\draw[dashed] (-1.4,1.5)--(-1.4,-0.5);
\draw[dashed] (1.4,1.5)--(1.4,-0.5);
\draw[<-] (-1.4, -0.3)--(-0.5,-0.3);
\draw[->] (0.5,-0.3)--(1.4, -0.3);
\node at (0,-0.3) {$\scriptstyle\blue{\ell}$};
\end{tikzpicture}
\end{center}
\caption{\footnotesize{{\bf Graph surgery, step 1:} A soliton is cut at the height
  such that the
  corresponding width equals the length of the pendant in the graph
  $\mathcal P$.}}
\end{figure}
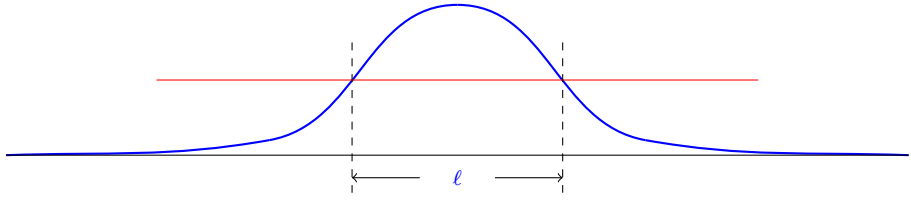

\begin{figure}[t] \label{surg1-bis}
\begin{center}
\begin{tikzpicture}[xscale=0.7,yscale=1]
\draw[-,thick,blue] (1,1.2) to [out=305,in=179] (5,0.02);
\draw[-,thick,blue] (-1,1.2) to [out=235,in=2] (-5,0.02);
\draw[-|] (-5,0) -- (-1,0);
\draw[|-] (1,0) -- (5,0);
\node at (-1.2,-0.4) {\blue{$\scriptscriptstyle -\ell/2$}};
\node at (1.2,-0.4) {\blue{$\scriptscriptstyle\ell/2$}};
\draw[-,thick,blue] (7,1) to [out=50,in=180] (8,1.8);
\draw[-,thick,blue] (8,1.8) to [out=0,in=130] (9,1);

\draw[|-|] (7,0) -- (9,0);
\node at (6.8,-0.4) {\blue{$\scriptscriptstyle -\ell/2$}};
\node at (9.2,-0.4) {\blue{$\scriptscriptstyle\ell/2$}};
\node at (14,0) {};
\end{tikzpicture}
\end{center}
\caption{\footnotesize{{\bf Pieces of soliton.} After the previous cut, one is left
  with three pieces: two tails and one head.}}
\end{figure}
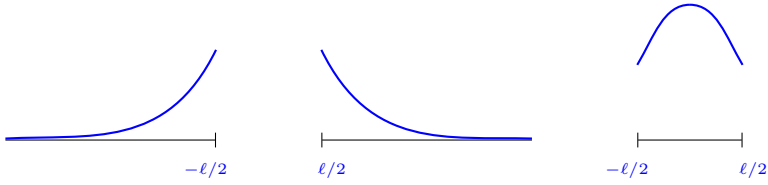

\noindent
Now we join the {two} resulting soliton {tails} (see
Fig.10) together at their maximum, and
  place them on the line of the graph $\mathcal P$, with the maximum at
  the vertex, as illustrated in Fig.11.

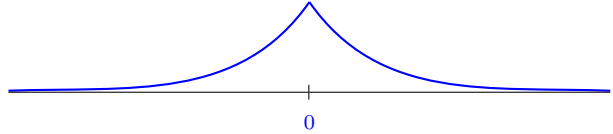
\begin{figure} \label{surg2}
\begin{center}
\begin{tikzpicture}[xscale=1,yscale=1]
\draw[-,thick,blue] (0,1.2) to [out=305,in=178] (4,0.02);
\draw[-] (0,0) -- (4,0);
\draw[-,thick,blue] (0,1.2) to [out=235,in=2] (-4,0.02);
\draw[-|] (-4,0) -- (-0,0);
\node at (0,-0.4) {$\blue{\scriptstyle 0}$};
\end{tikzpicture}
\end{center}
\caption{\footnotesize{{\bf Graph surgery, step 2.} The two tails produced by the
  first surgery step and shown in Fig. \ref{surg1-bis} are pasted
  together in such a way that the maximum is located at the
  vertex. Notice also that the maximum is a corner point.}}
\end{figure}

\noindent
Then, we {rearrange} the {head} of the soliton displayed in Fig. \ref{surg1-bis}
{ monotonically} on the pendant, namely on the interval {$[0, \ell]$}
(see Fig.12). Notice that for every point on the range,
except the maximum, the number of preimages passes from two to
one. This makes the energy strictly decrease.

\begin{figure} \label{rear}
 \begin{center}
\begin{tikzpicture}[xscale= 1,yscale=1]
\draw[|-|] (-2,0) -- (2,0);
\node at (-2,-0.4) {$\blue{\scriptstyle 0}$};
\node at (2,-0.4) {$\blue{\scriptstyle \ell}$};
\draw[-,thick,blue] (-2,1) to [out=30,in=180] (2,2);
\draw[-,thick,blue] (-9,1) to [out=50,in=180] (-7,2);
\draw[-,thick,blue] (-7,2) to [out=0,in=130] (-5,1);
\draw[|-|] (-9,0) -- (-5,0);
\node at (-9.2,-0.4) {\blue{$\scriptscriptstyle -\ell/2$}};
\node at (-5.2,-0.4) {\blue{$\scriptscriptstyle\ell/2$}};
\node at (-3.5,1) {\red{$\longrightarrow$}};
\end{tikzpicture}
\end{center}
\caption{\footnotesize{{\bf Monotone rearrangement on the pendant.} The monotone
  rearrangement of the soliton head lowers the energy.}}
\end{figure}
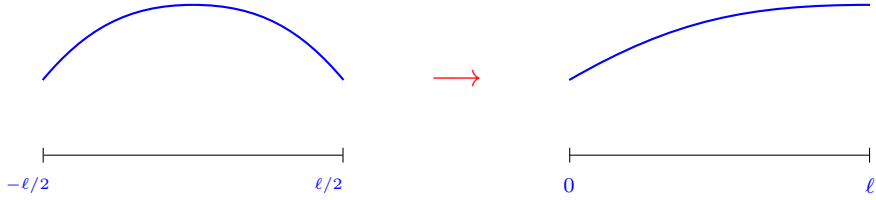
\noindent
Finally, we put all pieces together and construct a function
$\widetilde u$ on
$\mathcal P$ (Fig.13), such that
$$
E(\widetilde u,\G) < - \frac {\mu^3}{96}.
$$
Then, by Proposition \ref{comparison}, the infimum is {attained}.

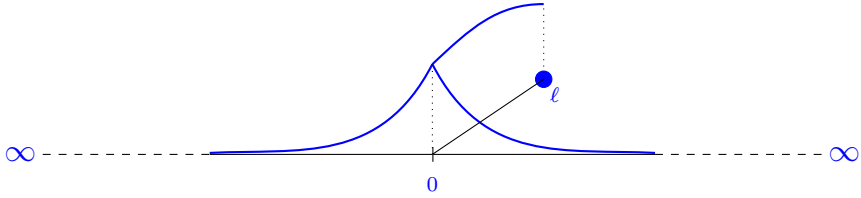
\begin{figure}
\begin{center}
\begin{tikzpicture}[xscale=0.74,yscale=1]
\node[nodo,blue] at (2,1){};
\node at (-7.4,0) [infinito]   {$\blue\infty$};
\node at (7.4,0) [infinito]   {$\blue\infty$};
\draw[-,thick,blue] (0,1.2) to [out=305,in=178] (4,0.02);
\draw[|-] (0,0) -- (4,0);
\draw[-,thick,blue] (0,1.2) to [out=235,in=2] (-4,0.02);
\draw[-] (-4,0) -- (-0,0);
\draw[-,thick,blue] (0,1.2) to [out=35,in=180] (2,2);
\node at (0,-0.4) {$\blue{\scriptstyle 0}$};
\node at (2.2,0.8){$\blue{\scriptstyle \ell}$};
\draw[dashed] (4,0) -- (7,0);
\draw[dashed] (-7,0) -- (-4,0);
\draw[-] (0,0) -- (2,1);
\draw[dotted] (2,1)--(2,2);
\draw[dotted] (0,0)--(0,1.2);
\end{tikzpicture}
\end{center}
\caption{\footnotesize{{\bf The competitor $\widetilde u$.} Through the procedure we illustrated a
  new function $\widetilde u$ is constructed, that turns out to be a
  good competitor for the soliton.}}
\end{figure}

\noindent
In this example, metric turns out to have no role in the existence of
the ground state. So, this
result relies on topology only.

\noindent
In the next section we  shall examine a case in which topology is not
sufficient to establish the existence of a ground state.













\noindent

\section{Metric} \label{sec:metric}
In the preceding section we stressed the role of topology in the
existence of ground states. But for graphs where assumption (H) fails,
existence and nonexistence can be conditioned by the metric. Let us
stress that,
according to the correspondence between long edges and large mass, the
example we are giving can be understood as a case in which, for a
fixed graph, the existence of a ground state depends on the choice of
the mass. We  will adopt the first point of view, in which the
existence depends on the metric, because it is more intuitive.

\noindent
Let $\mathcal G_\ell$ be a graph made of three halflines and a finite pendant of
length $\ell$ (see Fig.14). We will try
to follow the reasoning applied in the case of the line with the
pendant, and show why it fails.
\begin{figure}[h] \label{gielle}
\begin{center}
\begin{tikzpicture}[xscale= 1.6,yscale=1.8]
\node at (-2.5,0) [infinito]  (1) {$\infty$};
\node at (0,0) [nodo] (2) {};
\node at (0.3,0.3) {$\vv$};
\node at (2.5,0) [infinito]  (3) {$\infty$};
\node at (1.5,0)  [minimum size=0pt] (5) {};
\node at (-1.5,0) [minimum size=0pt] (6) {};
\node at (0,1) [nodo] (7) {};
\node at (0,-1.5) [infinito] (8) {$\infty$};
\node at (0,-0.8) [minimum size=0pt] (9) {};
\draw[dashed] (1) -- (6);
\draw[dashed] (3) -- (5);
\draw [-] (2) -- (6) ;
\draw [-] (5) -- (2) ;
\draw [-] (2) -- (7) ;
\draw [-] (2) -- (9) ;
\draw[dashed] (9) -- (8);
\end{tikzpicture}
\end{center}
\caption{\footnotesize{{\bf The graph $\mathcal G_\ell$:} three halflines and a
  finite pendant of length $\ell$ that meet together at a unique
  vertex $\vv$.}}
\end{figure}
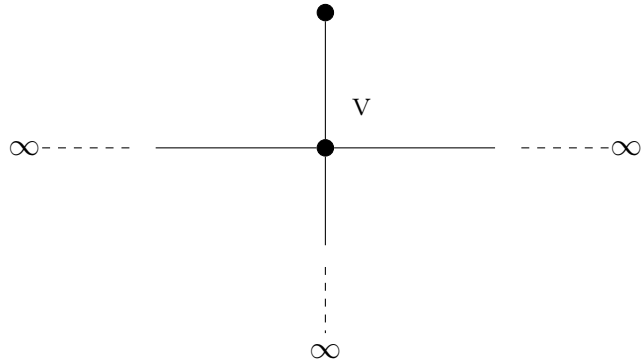
It immediately appears that, fixed $\mu$, if the
pendant is long enough, then it is possible to place almost exactly
half a soliton on it, producing then a function whose energy is less
than the energy of the soliton, so, by Proposition \ref{comparison},
the infimum is attained and then a ground state exists.

\noindent
Conversely, if $\ell$ is small enough, then one can prove that a
ground state does not exist: indeed, owing to Theorem 4.3 in
\cite{ast-arxiv}, if there is a ground state for every $\ell > 0$, then
there is a ground state for $\ell = 0$ too. But for $\ell = 0$ the
graph $\G$ reduces to the star graph made of three halflines, that
falls into the scope of Theorem \ref{nonexistence}. It follows that,
as the pendant grows, there is at least one transition from existence
and nonexistence. Let us show that there is {\em one} transition
only, i.e. that there exist a unique {\em critcal length $\ell^\star >
  0$} such that a ground state exists if and only if $\ell \geq \ell^\star$.

To this aim,
let $\ell$ be such that a ground state $\psi_\ell$ correspondingly
exists (Fig.15). Then, in analogy with the case of the line with a pendant, it
is possible to show, just arguing on the number of counterimages, that
the ground state looks like as shown in Fig.15.
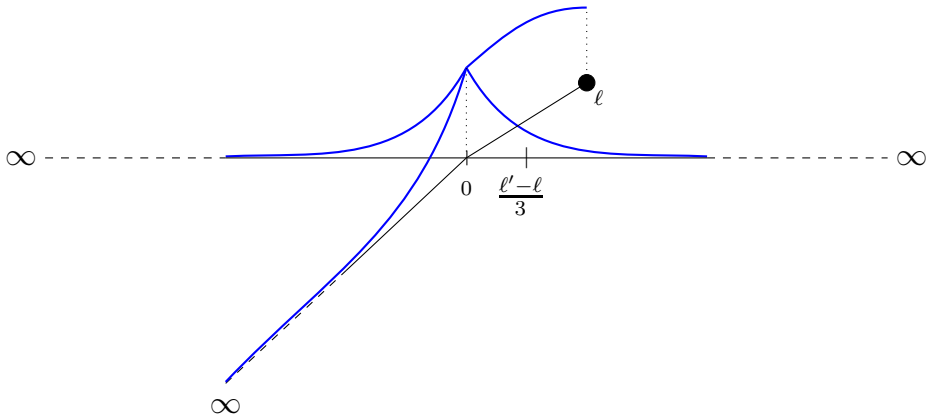
\begin{figure} \label{candidate}
\begin{center}
\begin{tikzpicture}[xscale= 0.8,yscale=1]
\node[nodo] at (2,1){};
\node at (-7.4,0) [infinito]   {$\infty$};
\node at (7.4,0) [infinito]   {$\infty$};
\node at (1,0) {$\scriptstyle{|}$};
\node at (0.9,-0.5) {$\frac{\ell' - \ell} 3$};
\node at (-4,-3.3) [infinito]   {$\infty$};
\draw[-][blue,thick](0,1.2) to [out=305,in=178] (4,0.02);
\draw[|-] (0,0) -- (4,0);
\draw[-][blue,thick] (0,1.2) to [out=235,in=2] (-4,0.02);
\draw[-] (-4,0) -- (-0,0);
\draw[-][blue,thick] (0,1.2) to [out=35,in=180] (2,2);
\node at (0,-0.4) {${\scriptstyle 0}$};
\node at (2.2,0.8){${\scriptstyle \ell}$};
\draw[dashed] (4,0) -- (7,0);
\draw[dashed] (-7,0) -- (-4,0);
\draw[dashed] (-4,-3) -- (-2,-1.5);
\draw[-] (-2,-1.5) -- (0,0);
\draw[-] (0,0) -- (2,1);
\draw[-,blue,thick](-4,-2.98) to [out=40,in=250] (0,1.2);
\draw[dotted] (2,1)--(2,2);
\draw[dotted] (0,0)--(0,1.2);
\end{tikzpicture}
\end{center}
\caption{\footnotesize{{\bf The ground state $\psi_\ell$.}
Assuming that the pendant is long enough to trap a ground state,
one can prove that the shape of such a ground state must be the one displayed here, with
the maximum at the tip of the pendant.
}}
\end{figure}
Consider now the graph ${\G_{\ell'}}$ with ${\ell' >
  \ell}$. Then, cut on every halfline the interval $[0, (\ell' -
  \ell)/3]$ and rearrange the three pieces monotonically together on
the interval {$[0, \ell'-\ell]$} (see Fig.16).

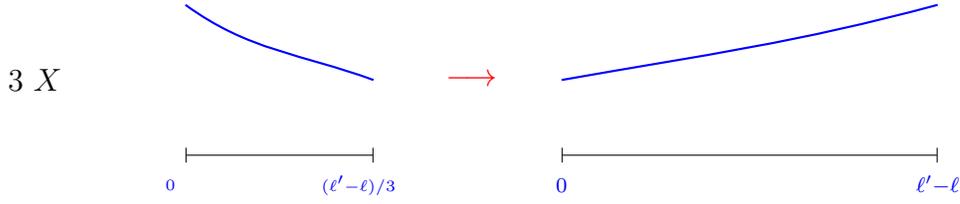
\begin{figure} \label{rear3}
\begin{center}
\begin{tikzpicture}[xscale=1,yscale=1]
\node at (-9,1) {$3 \ X $};
\draw[|-|] (-2,0) -- (3,0);
\node at (-2,-0.4) {$\blue{\scriptstyle 0}$};
\node at (3,-0.4) {$\blue{\scriptstyle \ell'-\ell}$};
\draw[-][blue,thick] (-2,1) to [out=10,in=195] (3,2);
\draw[-] [blue,thick] (-7,2) to [out=324,in=160] (-4.5,1);
\draw[|-|] (-7,0) -- (-4.5,0);
\node at (-7.2,-0.4) {\blue{$\scriptscriptstyle 0$}};
\node at (-4.7,-0.4) {\blue{$\scriptscriptstyle(\ell'-\ell) / 3$}};
\node at (-3.2,1) {\red{$\longrightarrow$}};
\end{tikzpicture}
\end{center}
\caption{\footnotesize{{\bf Rearranging the top of the tails of $\psi_\ell$.} On the exceeding
  length of the pendant of $\G_{\ell'}$ one puts the collective
  monotone rearrangement
of the tips of the three soliton tails of $\psi_\ell$, 
as portrayed in Fig.16. This step
lowers the energy level.}}
\end{figure}
Again, the loss of preimages yielded by the rearrangement lowers the
energy level, thus one can mount
the three pieces, obtaining a function whose energy is lower than the
energy of $\psi_\ell$.
\begin{figure}
\begin{center}
\begin{tikzpicture}[xscale= 0.8,yscale=1]
\node[nodo] at (2,1){};
\node[nodo] at (3,1.5){};
\node at (-7.4,0) [infinito]   {$\infty$};
\node at (7.4,0) [infinito]   {$\infty$};
\node at (-4,-3.3) [infinito]   {$\infty$};
\draw[-][blue,thick](0,1) to [out=305,in=178] (4,0.02);
\draw[|-] (0,0) -- (4,0);
\draw[-][blue,thick] (0,1) to [out=235,in=2] (-4,0.02);
\draw[-] (-4,0) -- (-0,0);
\draw[-][blue,thick] (0,1) to [out=15,in=200] (1,1.5);
\draw[-][blue,thick] (1,1.5) to [out=30,in=190](3,2.5);
\node at (0,-0.4) {${\scriptstyle 0}$};
\node at (3.1,1){${\scriptstyle \ell'}$};
\node at (1.2,0.3){${\scriptstyle \ell'-\ell}$};
\draw[dashed] (4,0) -- (7,0);
\draw[dotted] (3,1.5) -- (3,2.5);
\draw[dashed] (-7,0) -- (-4,0);
\draw[dashed] (-4,-3) -- (-2,-1.5);
\draw[-] (-2,-1.5) -- (0,0);
\draw[-] (0,0) -- (3,1.5);
\draw[-,blue,thick](-4,-2.98) to [out=40,in=250] (0,1);
\draw[dotted] (2,1)--(2,2);
\draw[dotted] (0,0)--(0,1.2);
\end{tikzpicture}
\end{center}
\caption{\footnotesize{{\bf A better competitor.} The function $\widetilde \psi$,
  resulting from the surgery and rearrangement procedure illustrated
  for $\psi_\ell$, has an energy content that is lower than the one of
  $\psi_\ell$. Again, the maximum of $\widetilde \psi$ is located at the tip of the pendant.}}
\end{figure}
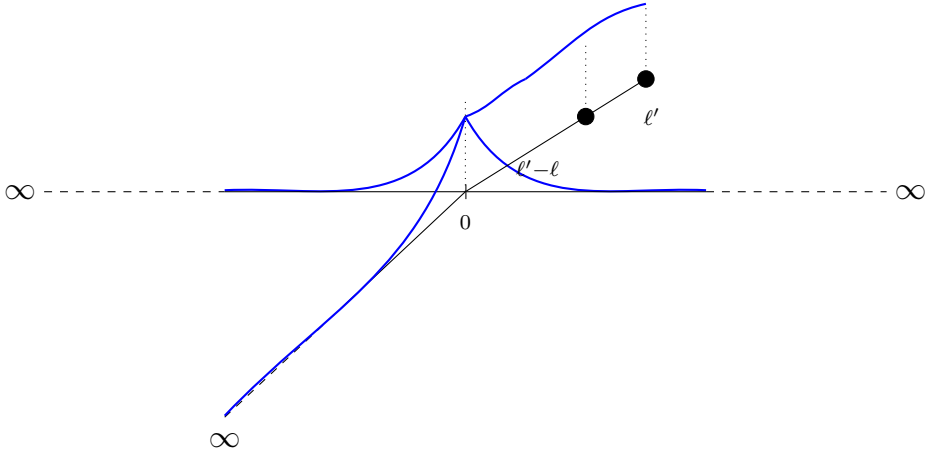

\noindent
In this way we obtained a function $\widetilde \psi \in H^1_\mu
(\G_{\ell'})$ (Fig.17) such that
$$
E ( \widetilde \psi, \G_{\ell'}) \ < \ E ( \psi_\ell, \G_{\ell}) \ \leq \ -
\frac {\mu^3}{96},
$$
and therefore a ground state at mass $\mu$ exists for the graph
$\G_{\ell'}$.

\noindent
Besides,
if $\ell' > \ell$, then
$$ \EE_{\G_{\ell'}} (\mu)  \ < \  \EE_{\G_{\ell}} (\mu),$$
so that $\EE_{\G_{\ell}} (\mu)$ is a monotonically decreasing function
of $\ell'$. Thus, once $\EE_{\G_{\ell}} (\mu) \leq - \mu^3 / 96$, as
$\ell$ increaes a
ground state still exists. Moreover, according to our analysis in the
case of large $\ell$, one has
$$
\EE_{\G_{\ell}} (\mu) \ \longrightarrow \ - \f {\mu^3}{24}, \qquad
\ell \to \infty.
$$

\section{Perspectives}
The problem of the existence of the ground state on graphs is physically relevant, for instance
in the context of Bose-Einstein condensates. Our results show that there is a competition between
the tendency of the condensate to stay on the compact core of the trap and the tendency to run away along 
long leads. The topology of the graph can solve between alternative: in particular,
the condensate can choose the first option only if in the graph hypothesis (H) is not satified,
that means that something like a bottleneck is present. 

\noindent
The analysis on the existence of ground states for the energy of a system living on a graph, can be extended in many
different directions: first of all, we remark once again that the same results
illustrated here for the cubic case hold for all subcritical
nonlinearity powers $\| u \|_{L^p (\G)}^p$ with $2 < p < 6$. However,
for non-cubic nonlinearities it is not possible to perform exact
computations, and then to write down the value of the energy level in such a precise form.

\noindent
A different behaviour is expected for the so-called {\em critical case}, that coincides with the 
Schr\"odinger equation with a quintic nonlinearity. This case is somewhat trivial
in $\R^N$ due to the presence of an additional symmetry, but we expect that the topology of the graph could somewhat
inhibit, at least partially, the new symmetry, leading to new and unexpected effects (this is
the main idea behind a work in progress with E. Serra and P. Tilli).

\noindent
Richer structures can be constructed by joining together pieces of different dimensions. On
the other hand, higher dimensions could be obtained by a limit procedure on the structure of graphs becoming more and more
"dense". Where should a ground state locate, in that case?

\noindent
Resarch on nonlinear effects on graphs is still at its beginning. Many developments are expected in 
the next future, both in the mathematics and phyiscs: for instance, in modelling periodic structures and/or cells, with the possible presence of a magnetic fields.

\bigskip

\noindent
{\bf Acknowledgements.}
 The author is partially supported by the
FIRB 2012 project ``Dispersive dynamics: Fourier Analysis and
Variational Methods'', code RBFR12MXPO, by the PRIN 2012 project
project ``Aspetti variazionali e perturbativi nei problemi
differenziali nonlineari'' and by the 2015 GNAMPA project ``Propriet\`a
spettrali delle equazioni di Schr\"odinger lineari e nonlineari''. 

The author is also grateful to Enrico Serra and Paolo Tilli for an uncountable number of
enlightening discussion and for many of the pictures of this paper.


\begin{thebibliography}{99}

\bibitem{acfn-rmp}
R.~Adami, C.~Cacciapuoti, D.~Finco,  D.~Noja, \emph{Fast solitons on star
  graphs}, Rev. Math. Phys \textbf{23} (2011), no.~4, 409--451.


\bibitem{acfn-jpa} R.~Adami, C.~Cacciapuoti, D.~Finco, D.~Noja,
\emph{On the structure of critical energy levels for the cubic focusing {NLS} on star graphs},
 J. Phys. A \textbf{45} (2012), no. 19, 192001, 7pp.

\bibitem{acfn-epl} R.~Adami, C.~Cacciapuoti, D.~Finco, D.~Noja,
{\em Stationary states of NLS on star graphs}, Europhys. Lett.
\textbf{100} (2012), no. 1,  10003, 6pp.


\bibitem{acfn-aihp} R. Adami, C. Cacciapuoti, D. Finco, D. Noja,
{\em Constrained energy minimization and orbital stability for the NLS equation on a star graph},  
Ann. Inst. Poincar\'e, An. Non Lin. \textbf{31} (2014), no. 6, 1289--1310. 

\bibitem{acfn-jde} R.~Adami, C.~Cacciapuoti, D.~Finco, D.~Noja, \emph{Variational properties and orbital stability of standing waves
  for {NLS} equation on a star graph},  J. Diff. Eq. \textbf{257} (2014), no. 10, 3738--3777.

\bibitem{an-jpa} R. Adami, D. Noja, 
{\em Existence of dynamics for a 1d NLS with a generalized point defect},
J. Phys. A: Math. Theor. \textbf{42} (2009), 495302, 19pp.



\bibitem{ast-cv}
R.~Adami, E.~Serra, P.~Tilli,
{\em NLS ground states on graphs},
Calc. Var. and PDEs,  to appear.  ArXiv: 1406.4036 (2014).

\bibitem{ast-arxiv}
R.~Adami, E.~Serra, P.~Tilli,
{\em Threshold phenomena and existence results for NLS ground states
  on graphs},
preprint arXiv:1505.03714 (2015).

\bibitem{ali}
F.~Ali Mehmeti.
{\em Nonlinear waves in networks.} Akademie Verlag, Berlin, 1994.


\bibitem{bellazzini}
B Bellazzini, M Mintchev, {\em Quantum fields on star graphs}
J. Phys. A: Math. Gen. 39 (35), 11101.

\bibitem{vonbelow}
J.~von~Below,
{\em An existence result for semilinear parabolic network equations with dynamical node conditions},
In Pitman Research Notes in Mathematical Series 266, Longman, Harlow Essex, 1992,  274--283.

\bibitem{berkolaiko}
G.~Berkolaiko, P.~Kuchment,
{\em Introduction to quantum graphs.}
Mathematical Surveys and Monographs, 186. AMS, Providence, RI, 2013.


\bibitem{exner}
J. Blank, P. Exner, M. Havlicek, {\em Hilbert space operators in Qunatum Physics}, Springer, New York, 2008.

\bibitem{bulgakov}
E. Bulgakov, A. Sadreev, {\em Symmetry-breaking in T-shaped photonic waveguides coupled with two identical nonlinear cavities}, Phys. Rev. B, 84:155304, 9pp, 2011.



\bibitem{cacciapuoti}
C.~Cacciapuoti, D.~Finco, D.~Noja,
{\em Topology induced bifurcations for the NLS on the tadpole graph},
Phys. Rev. E {\bf 91} (2015), no.1, 013206, 8pp. 


\bibitem{caudrelier}
V. Caudrelier, {\em On the Inverse Scattering Method for Integrable
  PDEs on a Star Graph}, Commun. Math. Phys. \textbf{338} (2015), no. 2,  893--917.


  
 \bibitem{friedlander}
 L. Friedlander, {\em Extremal properties of eigenvalues for a metric graph }, Ann. Inst. Fourier,
 55 (1), 199-211, 2005.
 
 



\bibitem{smi}
S.~Gnutzman, U.~Smilansky, S.~Derevyanko,
{\em Stationary scattering from a nonlinear network},
Phys. Rev. A {\bf 83} (2001), 033831, 6pp.



\bibitem{kevrekidis}
P.G.~Kevrekidis, D.J.~Frantzeskakis, G.~Theocharis, I.G.~Kevrekidis.
{\em Guidance of matter waves through Y-junctions},
Phys. Lett. A {\bf 317} (2003), 513--522.

\bibitem{kuchment}
P.~Kuchment,
{\em Quantum graphs I. Some basic structures},
Waves in Random Media {\bf 14} (2004), no. 1, S107--S128.




\bibitem{marzuola} J.L. Marzuola, D.E. Pelinovsky, 
{\em Ground state on the dumbbell graph}, arXiv:1509.04721 (2015).



\bibitem{matrasulov}
Z.~Sobirov, D.~Matrasulov, K.~Sabirov, S.~Sawada, K.~Nakamura.
\emph{Integrable nonlinear {S}chr\"{o}dinger equation on simple networks: connecion formula at vertices},
Phys. Rev. E \textbf{81} (2010), no. 6, 066602, 10pp.

\bibitem{noja14}
D.~Noja.
{\em Nonlinear Schr\"odinger equation on graphs: recent results and open problems},
Philos. Trans. R. Soc. Lond. Ser. A Math. Phys. Eng. Sci.  {\bf 372} (2014), no. 2007, 20130002, 20pp.

\bibitem{NPS} D.~Noja, D.~Pelinovsky, G.~Shaikhova, {\em Bifurcation and stability of standing waves in the nonlinear Schr\"odinger equation on the tadpole graph}, Nonlinearity, {\bf 28}, (2015) 2343-2378 


\bibitem{t16}
L. Tentarelli, {\em NLS ground states on metric graphs with localized nonlinearities},  J. Math. An. and Appl., 433 (1), 201--304, 2016.


\bibitem{tokuno} A. Tokuno, M. Oshikawa, E. Demler, {\em Dynamics of the one-dimensional Bose liquids: Andreev-like reflection at Y-junctions and the absence of Aharonov-Bohm effect}, Phys. Rev. Lett., 100:140402, 4pp,
2008.

\bibitem{hannes}
H. Uecker, D. Grieser, Z.~Sobirov, D. Babajanov, D.~Matrasulov, 
{\em Soliton transport in tubular networks: Transmission at vertices in the shrinking limit},
Phys. Rev. E {\bf 91} (2015), no. 2, 023209.

\bibitem{vidal}
E.J.G.~Vidal, R.P.~Lima, M.L.~Lyra.
{\em Bose-Einstein condensation in the infinitely ramified star and wheel graphs},
Phys. Rev. E  {\bf 83}  (2011), 061137, 8pp.


\bibitem{malomed}
N. Viet Hung, M. Trippenbach, B. Malomed, {\em Symmetric and asymmetric solitons trapped in H-shaped potentials},
Phys. Rev. A, 84:053618, 10pp. 2011.

\bibitem{zapata} I. Zapata, F. Sols, {\em Andreev reflection in bosonic condensates}, Phys. Rev. Lett., 102:180405. 4pp. 2009.


\end{thebibliography}
\end{document}